%

\magnification=\magstep1
\overfullrule=0pt
\input amstex
\documentstyle{amsppt}

\hsize=5.5truein
\vsize=8.5truein
\hoffset0.5truein

\define\CUB{\operatorname{CUB}}

\define\ORD{\operatorname{ORD}}

\define\lon{$L[O^\#]$}
\define\lan{\langle}
\define\ran{\rangle}

\document

\baselineskip=.25truein

\centerline{\bf Generic Saturation}
\vskip10pt

{\centerline{ Sy D. Friedman \footnote"*"{Research supported by NSF \#
9205530 DMS et l'Universite' de Marne la Vallee.}}} 
\centerline {M.I.T. et l'Universite' de Paris 7}
\centerline {sdf$\@$math.mit.edu}

\vskip20pt

Assuming that ORD is $\omega +\omega $-Erd\"os we show
that if a class forcing amenable to $L$ (an $L$-forcing) has a generic
then it has one definable in a set-generic extension of $L[O^\#]$.  In fact
we may choose such a generic to be {\it periodic} in the sense that
it preserve the indiscernibility of a final segment of a periodic
subclass of the Silver indiscernibles, and therefore to be {\it almost
codable} in the sense that it is definable from a real which is generic for
an $L$-forcing (and which belongs to a set-generic extension of
$L[O^\#]$). This result is best possible in the sense that for any
countable ordinal $\alpha $ there is an $L$-forcing which has generics but
none periodic of period $\le\alpha .$  However, we do not know if an
assumption beyond $ZFC+O^\#$ exists is actually necessary for these
results.

Let $P$ denote a class forcing definable over an amenable ground model
$\lan L,A\ran$  and assume that $O^\#$ exists.

\vskip10pt

\noindent
{\bf Definition.} \  $P$ is {\it relevant} if $P$ has a generic definable
in \lon. $P$ is {\it almost relevant} if $P$ has a generic definable in a
set-generic extension of \lon.

\vskip10pt

\noindent
{\bf Remark.} \ The reverse Easton product of Cohen forcings $2^{<\kappa
}$. $\kappa $ regular is relevant. So are the Easton product and the full
product, provided $\kappa $ is restricted to the successor cardinals. See
Chapter 3, Section Two of Friedman [97]. Of course any set-forcing (in $L)$
is almost relevant.

\vskip10pt

\noindent
{\bf Definition.} \ $\kappa $ is $\alpha $-Erd\"os if whenever $C$ is CUB
in $\kappa $ and $f:[C]^{<\omega }\longrightarrow \kappa $ is regressive
(i.e., $f(a)<\min (a))$ then $f$ has a homogeneous set of ordertype
$\alpha.$

\vskip10pt

\noindent
{\bf Definition.} \ Let $\Cal{A}=\lan T,\epsilon ,\dots\ran$ be transitive
(in a countable language). $I\subseteq \ORD(T)$ is a {\it good set of 
$\Sigma_1$ indiscernibles} for $\Cal{A}$ if $\gamma \in I\longrightarrow I-\gamma $ is
a set of $\Sigma_1$ indiscernibles for $\lan \Cal{A},\alpha \ran_{\alpha <\gamma }.$

\vskip10pt

\noindent
{\bf Fact.} \ $\kappa$ is $\alpha $-Erd\"os iff whenever $\Cal{A}=\lan
T,\epsilon ,\dots\ran$ is transitive (in a countable language), $\kappa
\subseteq \ORD(T), C \CUB$ in $\kappa $ then there exists $I\subseteq
C,$ ordertype $(I)=\alpha $ such that $I$ is a good set of $\Sigma_1$ indiscernibles
for $\Cal{A}$.

\vskip10pt

\noindent
{\bf Theorem 1.} \ {\it Suppose $P$, defined over $\lan L,A\ran$, has a generic $G$ and there is a good
set $X$ of $\Sigma_1$ indiscernibles for $\lan L[O^\#,G], \epsilon ,G,A\ran$ of ordertype $\omega +\omega $ such that $\alpha \in X \longrightarrow \alpha$ 
is $\Sigma_1-stable$ in $O^\#,G,A$ (i.e., $\lan L_\alpha[O^\#,G],\epsilon,
G\cap L_\alpha,A\cap L_\alpha\ran$ is $\Sigma_1-$ elementary in $\lan L[O^\#,G],\epsilon,G,A\ran$). Then $P$ is almost relevant.}

\vskip10pt

\noindent
{\bf Corollary 2.} \ {\it Suppose $P$ has a generic and ORD is $\omega
+\omega $-Erd\"os. Then $P$ is almost relevant.}

\vskip10pt
\noindent
{\bf Remark.} If $\{\kappa|\kappa$ is $\alpha$-Erd\"os$\}$ is stationary then it 
follows that ORD is $\alpha$-Erd\"os.
\vskip10pt
The proof of Theorem 1 provides a stronger conclusion which we describe next.

\vskip10pt

\noindent
{\bf Definition.} \ $P$ is {\it codable} if $P$ has a generic $G$ definable
over $L[R],$ $R$ a real generic over $L,$ $R\in$ \lon. $P$ is {\it almost}
codable if $P$ has a generic $G$ definable over $L[R],R$ a real generic
over $L,R$ in a set-generic extension of \lon.

These notions can be alternatively described in terms of
indiscernibility-preservation:

\vskip10pt

\noindent
{\bf Definition.} \  Let $I=\lan i_\gamma | \gamma
\in \ORD\ran$ be the increasing enumeration of the Silver
indiscernibles. For any ordinals $\lambda _0,\lambda (\lambda >0)$ define
$I_{\lambda _0,\lambda }=\{i_\alpha |\alpha $ of the form $\lambda
_0+\lambda \cdot\beta ,\beta \in\ORD\}.$ $P$ is $\lambda _0,\lambda $-{\it
periodic} if there is a $P$-generic $G$ such that $I_{\lambda _0,\lambda }$
is a class of indiscernibles for $\lan L[G],\epsilon ,G,A\ran.$  $P$ is
{\it almost} $\lambda _0,\lambda $-{\it periodic} if it is $\lambda
_0,\lambda $-periodic in a set-generic extension of $V.$

\vskip10pt

\noindent
{\bf Proposition 3.} \

(a) \ If $A=\emptyset, P$ $L$-definable without parameters then $P$ is
codable iff $P$ is almost $\lambda _0,\lambda $-periodic for some $\lambda _0,\lambda .$

(b) \ $P$ is almost codable iff $P$ is almost $\lambda _0,\lambda
$-periodic for some $\lambda _0,\lambda.$

\vskip10pt

\flushpar
{\smc Proof.} \  (a) \ For the ``only if'' direction, see Chapter 5,
Section Two of Friedman [97]. For the ``if'' direction, build a tree in
\lon, a branch through which produces a real coding a generic witnessing
$\lambda _0,\lambda $-periodicity for some countable $\lambda _0,\lambda $.
Then this tree has a branch in \lon, proving that $P$ is codable.  Part (b)
is similar (and does not need the assumption of (a) since for any
$A$ there exists some $\lambda _0$ such that $I_{\lambda _o,1}$ is a class
of indiscernibles for $\lan L,A\ran).$  \hfill{$\dashv$}

\vskip10pt

\noindent
{\bf Remark.} It follows that in Theorem 1 and Corollary 2, if $A=\emptyset,
P$ is $L-$definable without parameters then ``almost relevant'' can be
replaced by ``relevant''.

\vskip10pt

The standard examples of relevant class forcing are in fact $0,1$-periodic.

\vskip10pt

\noindent
{\bf Periodicity Conjecture.} \ If $P$ has a generic then $P$ is almost
$\lambda _0,\lambda $-periodic for some countable $\lambda .$

Our proof of Theorem 1 establishes the Periodicity Conjecture, under the
extra hypothesis that ORD is $\omega + \omega$-Erd\"os:

\vskip10pt

\noindent
{\bf Theorem 4.} \ {\it Suppose $P$ satisfies the hypothesis of Theorem 1.
Then P is almost $\lambda _0,\lambda $-periodic for some
countable $\lambda .$ Thus if ORD is $\omega+\omega$-Erd\"os then the 
Periodicity Conjecture is true.}

The Periodicity Conjecture cannot be strengthened.

\vskip10pt

\noindent
{\bf Theorem 5.} \ {\it Suppose $\alpha ,\beta $ are ordinals, $\beta $
countable. Then there is an $L$-forcing $P$ such that P has a generic but P
is not almost $\lambda _0,\lambda $-periodic for $\lambda _0<\alpha $ or
for $\lambda <\beta .$ }

\vskip10pt

\noindent
{\smc Proof of Theorem 4.} \ Fix a $P$-generic $G$ as in the hypothesis of
Theorem 4; we shall construct another $P$-generic $G^*$ such that for some
$\lambda _0$ and countable $\lambda ,$ $I_{\lambda _0,\lambda }$ is a class
of indiscernibles for $\lan L[G^*],\epsilon ,G^*,A\ran.$ Let $X$ be a good
set of $\Sigma_1$ indiscernibles for $\lan L[O^\#,G],\epsilon ,G,A\ran$ of ordertype
$\omega+\omega$ such that $\alpha\in X \longrightarrow \alpha$ is $\Sigma_1-$Stable in $O^\#,G,A$.

Select a canonical enumeration of the $\lan L,A\ran$-definable open dense
subclasses of $P:$ Thus let $\lan D_n|n\in\omega \ran$ be a sequence of
predicates where each $D_n(x,\alpha _1\dots\alpha _n)$ is definable over
$\lan L,A\ran$ such that for each $\alpha_1 <\dots<\alpha _n$ in ORD,
$\{x\in L|D_n(x,\alpha _1\dots\alpha _n)\}$ is an open dense subclass of
$P$ and every open dense subclass of $P$ is of this form for some $n$, for some $\alpha _1<\dots<\alpha _n$ in $I=$(Silver) indiscernibles. 
We may also assume that $\{\lan n,x,\vec\alpha \ran|D_n(x,\vec\alpha )\}$
is definable over $\lan L,A\ran$ relative to a satisfaction predicate for
$\lan L,A\ran.$ For $\alpha _1<\dots<\alpha _n$ in ORD we abuse notation
and write $D(\alpha_1\dots\alpha _n)$ for $\{x\in L|D_n(x,\alpha
_1,\dots\alpha _n)\}.$ Also let $D^*(\alpha_1\dots\alpha_n)=\cap\{D(\vec\beta)|\vec\beta\subseteq\vec\alpha\}$.

Now we construct an $\omega $-sequence of terms with indiscernible
parameters which we will use to define $G^*.$

For $j_0\in X$ choose the least $t_{j_{0}}(\vec k_0(j_0),j_0,\vec
k_1(j_0))$ in $D(j_0)\cap G,$ where $t_{j_{0}}$ is a Skolem term for
$L,\vec k_0(j_0)<j_0<\vec k_1(j_0)$ is an increasing sequence of
indiscernibles. By the good-indiscernibility of $X$, $t_{j_{0}}=t_0$, $\vec
k_0(j_0)=\vec k_0$ are fixed. Thus we can write $t_0(\vec k_0,j_0, \vec
k_1(j_0))\in D(j_0)\cap G$ for $j_0\in X.$ By the $\Sigma_1-$stability
in $O^\#,G,A$ of the elements of $X$ we have: $j_0<j_1$ in $X \longrightarrow
\vec k_1(j_0)<j_1.$

Next for $j_0<j_1$ in $X$ choose the least $t_{j_0,j_1}(\vec
k^1_0(j_0,j_1),j_0, \vec k^1_1(j_0,j_1), j_1, \vec k^1_2(j_0,j_1))$ in
$D^*(\vec k_0, j_0, \vec k_1(j_0), j_1, \vec k_1(j_1))\cap G.$ By the
good-indiscernibility of $X$ we can write the above term with indiscernible
parameters as $t_1(\vec k^1_0, j_0, \vec k^1_1(j_0), j_1, \vec k^1_2(j_0,
j_1)).$  However, we want to argue that $\vec k^1_2(j_0,j_1)$ can be chosen
independently of $j_0.$  To arrange this, first note that $t_{j_0,j_1}(\vec
k^1_0(j_0,j_1),j_0, \vec k^1_1(j_0,j_1), j_1, \vec k^1_2(j_0,j_1))$ 
= $t_{j_0,j_1}(\vec
k^1_0(j_0,j_1),j_0, \vec k^1_1(j_0,j_1),j_1,\vec
k_{2,0}^1(j_0,j_1),\vec\infty)$ where the latter is independent of the choice of the
indiscernibles $\vec\infty$ above $\vec k^1_{2,0}(j_0,j_1)$ and where
$(\vec k^1_0(j_0,j_1), \vec k^1_1(j_0,j_1), \vec k^1_{2,0}(j_0,j_1))$ is
the least sequence of {\it ordinals} such that this term with parameters
belongs to $D^*(\vec k_0,j_0,\vec k_1(j_0), j_1, \vec k_1(j_1))\cap G\cap
L_{\min \vec\infty}.$ By the
good-indiscernibility of $X$ we can write this as $t_1(\vec k^1_0, j_0,
\vec k^1_1(j_0),j_1,\vec k^1_{2,0}(j_0,j_1),\vec\infty).$ Note that
$(\vec k^1_0,\vec k^1_1(j_0), \vec k^1_{2,0}(j_0,j_1))$ is definable in
$\lan L[G], G,A\ran$ from $\vec\infty,\vec k_0,j_0$,$\vec k_1(j_0),j_1,\vec
k_1(j_1)$ and therefore $\vec k^1_{2,0}(j_0,j_1)$ is definable in $\lan
L[G], G,A\ran$ from $\vec\infty, \vec k_1(j_1)$ and parameters $\le j_1.$

\vskip10pt

\noindent
{\bf Lemma 6.} \ $\vec k^1_{2,0}(j_0,j_1)$ is independent of $j_0.$

\vskip10pt

\noindent
{\smc Proof.} \ Enumerate the first $\omega +1$ elements of $X$ in
increasing order as $j_0<j_1<\dots < j=(\omega +1)$st element of $X$ and
for any $m,n$ let $\vec k(j_n,j)$ $(m)$ denote the $m^{\text{th}}$ element
of $\vec k^1_{2,0}(j_n,j)$. If the Lemma fails then for some fixed $m,\vec
k(j_0,j)(m)<\vec k(j_1,j)(m)<\dots$ forms an increasing $\omega
$-sequence of indiscernibles with supremum $\ell\in I.$ By the remark
immediately preceding this Lemma, $\ell$ has cofinality $\le j$ in
$L[G]$. By Covering between $L$ and $L[G],\ell$ has cofinality $<(j^+$ in
$L[G])$ in $L$. This contradicts the following.

\vskip10pt

\noindent
{\bf Claim.} \ $j^+$ in $L[G]=j^+$ in $L$.

\vskip10pt

\noindent
{\smc Proof of Claim.} \ If not then in $L[G]$ there is a CUB $C\subseteq
j$ such that $C$ is almost contained in each constructible $D\subseteq j.$
But $I\cap j$ is the intersection of countably many such $D$ and therefore
as $j$ is regular (in $L[G,O^\#])$ we get that $C$ is
almost contained in $I;$ so $O^\#$ belongs to $L[G],$ contradiction. This
proves the Claim and hence the Lemma. \hfill{$\dashv$}

\vskip10pt

Thus we can write $t_1(\vec k^1_0, j_0,\vec k^1_1(j_0), j_1, \vec
k_2^1(j_1))\in $ $D^*(\vec k_0, j_0, \vec k_1(j_0), j_1, \vec k_1(j_1))\cap
G$ for $j_0<j_1$ in $X$. By modifying the term $t_1$ we may assume that
$\vec k^1_1(j_0)=\vec k^1_2(j_0)$ for $j_0\neq \min (X).$ Also we can
assume that $\vec k_0\subseteq \vec k^1_0, \vec k_1(j_0)\subseteq \vec
k^1_1(j_0)$ for $j_0\in X$ and moreover that the structure $\lan \vec
k_1^1(j_0),<\ran$ with a unary predicate for $\vec k_1(j_0)$ has
isomorphism type independent of $j_0\in X.$

We obtain $t_2$ in a similar way: \ thus,
$$\aligned
&t_2(\vec k^2_0, j_0, \vec k^2_1(j_0), j_1, \vec k^2_1(j_1), j_2, \vec
k^2_1(j_2))\in\\ 
&D^*(\vec k^1_0, j_0, \vec k^1_1(j_0), j_1, \vec k^1_1(j_1), j_2, \vec
k^1_1(j_2))\cap G\endaligned$$
for $j_0<j_1$ in $X$ and $\vec k^1_0\subseteq \vec k^2_0, \vec
k^1_1(j_0)\subseteq \vec k^2_1(j_0)$, $\lan \vec k^2_1(j_0), <\ran$ with
unary predicates for $\vec k^1_1(j_0),\vec k_1(j_0)$ has isomorphism type
independent of $j_0$. Continue in this way to define $t_n(\vec k^n_0, j_0, \vec k^n_1(j_0),\dots ,
j_n, \vec k^n_1(j_n))$ for each $n$ and for $j_0<\dots <j_n$ in $X.$ (The
analogous version of Lemma 6 uses the first $\omega+n$ elements of $X$.)

Let $i_{\lambda _{0}}=\min X$ and
$\lambda= $ ordertype $(\bigcup\limits_{n}\vec k^n_1(j_0))$ for $j_0\in X,$ an ordinal independent of the choice of $j_0.$

We may assume that $\lambda $ is a limit ordinal and in a generic extension
where $\lambda _0$ is countable we may arrange that $\bigcup\limits_{n}\vec
k^n_0=I\cap i_{\lambda _0}.$  Also note that $I-i_{\lambda _0}$ is a class
of indiscernibles for $\lan L,A\ran.$ Now in $V[g]$, where $g$ is a L\'evy
collapse of $i_{\lambda _0}$ to $\omega $ carry out the above construction,
arranging that $\bigcup\limits_{n}\vec k^n_0=i_{\lambda _0}.$ For any
indiscernible $i_\delta$ define $\vec k^n_1(i_\delta )\subseteq I\cap
(i_\delta , i_{\delta +\lambda })$ so that $\lan I\cap (i_\delta ,
i_{\delta +\lambda }),<\ran$ with a predicate for $\vec k^n_1(i_\delta )$
is isomorphic to $\left\lan \bigcup\limits_n\vec k^n_1(j_0),<\right\ran$
with a predicate for $\vec k^n_1(j_0),$ for $i_{\lambda _0}<j_0\in X.$ Define:
$$\aligned
G^*=\{p\in P|p&{\text { is extended by some }} t_n(\vec k^n_0, i_{\lambda
_1}, \vec k^n_1(i_{\lambda _1}),\dots\\
&i_{\lambda _n}, \vec k^n_1(i_{\lambda _{n}})) \text{ where } \lambda _0\le
\lambda _1<\dots<\lambda _n\\
&\text{ are of the form } \lambda _0+\lambda \cdot\alpha ,\alpha \in\ORD\}.
\endaligned$$
Using the indisciernibility of $I-i_{\lambda _0}$ in $\lan L,A\ran$ we see
that $G^*$ is compatible and meets every $\lan L,A\ran$-definable open
dense class on $P.$ Thus $P$ is $\lambda _0,\lambda $-periodic in $V[g].$
Note that $\lambda $ is countable in $V.$ This proves Theorem
4. \hfill{$\dashv$}

\vskip10pt

\noindent
{\bf Remark.} \ The proof of Theorem 4 only made use of a weaker
hypothesis: \ Define $X$ to be a good set  of {\it $\Sigma_1$ n-indiscernibles} for
$\Cal{A}=\lan T,\epsilon ,\dots\ran$ if $\gamma \in X\longrightarrow
X-\gamma $ is $\Sigma_1$ indiscernible for $\Cal{A}$ {\it for n-tuples}. Our proof
only used the existence of $X_1\supseteq  X_2\supseteq \dots$ such that each $X_n$ is a
good set of $\Sigma_1$ $n$-indiscernibles for $\lan L[O^\#,G], \epsilon , G, A\ran$ of
ordertype at least $\omega +\omega$ such that $\alpha \in X_n \longrightarrow \alpha$ is $\Sigma_1-$stable in $O^\#,G,A.$ This hypothesis is weaker in terms of
consistency strength than the hypothesis stated in Theorem 4.

\vskip10pt

\noindent
{\smc Proof of Theorem 5.} \ We employ here the techniques of Friedman [90]
and Friedman [94]. In the former, an $L$-definable forcing is constructed
so as to have a unique generic, which can be considered to be a real. In
Friedman [97], Chapter 5, Section Two it is shown that there exist reals
$R$ such that $I^R=$ Silver indiscernibles for $L[R]$ is equal to Even
$(I)=\{i_{2\alpha }|\alpha \in\ORD\}.$ By combining the latter construction
with the construction of Friedman [90] one obtains an $L$-definable forcing
$Q$ with a unique generic real $R,$ such that $I^R=$ Even $(I)$.

Now suppose that $\alpha $ is an $L$-countable ordinal. Define an iterated
class forcing as follows: $P_0=\{0\}.$ $P_{\beta +1}=P_\beta *P(\beta )$
where $P(\beta )$ applies the forcing $Q^{R_\beta }=(Q$ relativized to
$R_\beta )$ over the model $L[R_\beta ],$ where $R_\beta =$ the $P_\beta
$-generic real. (Thus if $R_{\beta +1}=$ the $P_{\beta +1}$-generic real we
get $I^{R_\beta ,R_{\beta +1}}=$ Even $(I^{R_\beta }).)$
For limit $\lambda \le\alpha $ let $P_\lambda =$ Inverse limit $\lan
P_\beta |\beta <\lambda \ran$ and $R_\lambda =$ Join of $\lan R_\beta
|\beta <\lambda \ran$ using the $L$-least counting of $\lambda .$

By Friedman [94], the $P_\beta$'s preserve cofinalities and ZFC. And
$P_\beta $-generics exist, using the methods of Friedman [97], Chapter 3,
Section Two. The forcing $P_\alpha $ adds a real $R$ such that
$I^R=\{i_{2\alpha \gamma }|\gamma \in \ORD\}$ (and has a unique generic).
If $\alpha $ is not countable in $L,$ first apply a L\'evy collapse of
$\alpha $ and then perform the above construction to obtain $P_\alpha .$
The generic is no longer unique (as the L\'evy collapse is not) but it is
the case that for any generic real $R,$ $I^R=\{i_{2\alpha \gamma }|\gamma
\in\ORD\}-(\alpha +1).$ To prove Theorem 5: Choose $\alpha $ to be the
$\beta $ of the statement of that theorem; then a $P_\alpha $-generic
exists (as $\alpha $ is countable) and $P_\alpha $ is not almost $\lambda
_0,\lambda$-periodic for $\lambda <\beta .$ To rule out the case $\lambda
_0<\alpha ^*=(\alpha $ of Theorem 5), add a Cohen set to $(\alpha ^*)^+$ of
$L$, after forcing with $P_\alpha.$ This proves Theorem 5. \hfill{$\dashv$}

\vskip10pt

\noindent
{\bf Questions.} \ (a) \ Is the Periodicity Conjecture provable in the
theory ZFC $+ O^\#$ exists?

(b) \ Suppose that whenever $P$ is an $L$-forcing with a generic $G$ such
that $\lan V[G], G\ran \vDash$ ZFC then there is such a $G$ definable in a
set-generic extension of $V.$  Does $O^{\#}$ exist?

(c) \ For which $\alpha $ countable in $L[O^\#]$ does there exist an
$L$-forcing $P$ with a unique generic $G,$ such that $\alpha $ is countable
in $L[G]?$

\vskip20pt

\centerline{\bf References}

\vskip10pt

Friedman [90] \ {\it The $\Pi^1_2$-Singleton Conjecture,} Journal of the
AMS, Volume {\bf 3,} Number 4.

\vskip5pt

Friedman [94] \ {\it Iterated Class Forcing,} Mathematical Research Letters
1.

\vskip5pt

Friedman [97] \ {\it Fine Structure and Class Forcing,} book manuscript.
\end{document}